\documentclass[11pt]{amsart}

\usepackage{amsmath,amssymb,amsthm,amscd,amsfonts,verbatim}
\usepackage[top=1in,bottom=1in,left=0.85in,right=0.85in]{geometry}
\usepackage{thmtools,thm-restate}
\usepackage{calc}
\usepackage{hyperref}

\title{Bounding the homology of FI-modules}

\author{Thomas Church}

\date{Original version December 28, 2015. References updated December 18, 2016.}

\theoremstyle{plain}

\declaretheorem[name=Theorem]{maintheorem}

\newcommand{\nc}{\newcommand}

\newcommand{\beq}{\begin{displaymath}}
\newcommand{\eeq}{\end{displaymath}}
\newcommand{\tensor} {\otimes}

\nc{\dmo}{\DeclareMathOperator}

\nc{\I}{\mathcal{I}}
\nc{\K}{\mathcal{K}}
\nc{\U}{\mathcal{U}}
\renewcommand{\L}{\mathbf{L}}
\nc{\ra}{\to}
\nc{\Q}{\mathbb{Q}}
\nc{\R}{\mathbb{R}}
\nc{\Z}{\mathbb{Z}}
\nc{\N}{\mathbb{N}}
\nc{\F}{\mathbb{F}}
\nc{\cF}{\mathcal{F}}
\nc{\A}{\mathbb{A}}
\nc{\T}{\mathcal{T}}
\nc{\tT}{\widetilde{\T}}
\nc{\V}{\mathcal{V}}
\nc{\PP}{\mathbf{P}}
\nc{\LL}{\mathbf{L}}
\nc{\cL}{\mathcal{L}}
\nc{\G}{\mathbb{G}}
\nc{\C}{\mathcal{C}}
\nc{\CP}{\C\PP}
\nc{\BB}{\mathcal{B}}
\nc{\sss}{s'}
\nc{\ttt}{t'}
\nc{\FF}{\mathcal{F}}
\dmo{\GL}{GL}
\dmo{\PSL}{PSL}
\dmo{\Teich}{Teich}
\dmo{\rank}{rank}
\nc{\gin}{i}
\nc{\ga}{\Gamma}
\dmo{\Conf}{Conf}
\dmo{\sign}{sign}
\dmo{\mult}{mult}
\dmo{\Fix}{Fix}
\dmo\Tr{Tr}
\dmo\Hom{Hom}

\dmo\coker{coker}

\dmo\im{im}
\dmo{\ima}{im}
\dmo\id{id}
\dmo\SL{SL}
\dmo\Sp{Sp}
\dmo\Mod{Mod}
\dmo\PMod{PMod}
\dmo\fd{fd}
\dmo\IA{IA}
\dmo\Schur{{\Bbb S}}
\dmo\Sym{Sym}
\dmo\Ind{Ind}
\dmo\Res{Res}
\dmo\tr{tr}
\dmo\Irr{Irr}
\dmo\st{st}
\dmo\Stab{Stab}
\dmo\End{End}

\nc{\bwedge}{{\textstyle\bigwedge}}
\nc{\case}[1]{\left\{\begin{array}{c}#1\end{array}\right.}
\nc{\abs}[1]{\left\lvert#1\right\rvert}

\dmo\FI{FI}
\dmo\FIMod{FI-Mod}
\dmo\FB{FB}
\dmo\FBMod{FB-Mod}
\dmo\dMod{-Mod}
\nc{\ZMod}{\Z\dMod}
\nc{\RMod}{R\dMod}
\nc{\FIsharp}{\FI\sharp}
\nc{\FIsharpMod}{\FIsharp\dMod}
\nc{\op}{\text{op}}
\nc\FIop{\FI^{\op}}
\nc\FIopMod{\FI^{\op}\dMod}
\dmo\Ob{Ob}
\dmo\Inj{Inj}
\dmo\Bij{Bij}
\dmo\Tor{Tor}

\renewcommand{\epsilon}{\varepsilon}
\nc{\coloneq}{\mathrel{\mathop:}\mkern-1.2mu=}
\nc{\margin}[1]{\marginpar{\scriptsize #1}}
\nc{\para}[1]{\medskip\noindent\textbf{#1.}}

\nc{\into}{\hookrightarrow}
\nc{\onto}{\twoheadrightarrow}

\nc{\comp}{\overline}
\nc{\cattag}{$\ast\ast$}
\nc{\cat}{\ref{catlabel}}

\nc{\arXiv}[1]{\href{http://arxiv.org/abs/#1}{arXiv:#1}}
\nc{\myemail}[1]{\href{mailto:#1}{\nolinkurl{#1}}}

\renewcommand{\phi}{\varphi}
\nc{\disjoint}{\sqcup}
\nc{\iso}{\simeq} 

\begin{document}
\begin{abstract}
The main theorem of Church--Ellenberg~\cite{CE} is a sharp bound on the homology of FI-modules, showing that the Castelnuovo--Mumford regularity of FI-modules over $\Z$ can be bounded in terms of generators and relations.
We give a new proof of this  theorem, inspired by earlier proofs by Li--Yu and Li. This should be read after \S2~of~\cite{CE}, which contains all needed terminology~and~background.
\end{abstract}
\maketitle
\vskip-10pt

The right-exact functors $H_0\circ D$ and $S\circ H_0$ from $\FIMod$ to $\FBMod$ are naturally isomorphic.\footnote{Being right-exact, it suffices to check that they agree on free/projective FI-modules, which on objects is Lemma 4.4 of \cite{CE}. On morphisms there is a bit more to check (since a map $M(X)\to M(Y)$ need not be induced by a map $X\to Y$), but it is 
pleasant to work this out by hand (our interest is instead in the map $X\to Y$ that it \emph{induces},~which~always~exists).}  Since $S$ is exact, 
\[S\circ \L H_0=\L (S\circ H_0)=\L (H_0\circ D)=(\L H_0)\circ (\L D).\]
Since $\L_q D=0$ for $q\geq 2$ \cite[Lemma 4.7]{CE}, 
we obtain for any $W$ a two-row spectral sequence converging to $S H_\ast W$.
Setting $K\coloneq \L_1 D(W)$, this two-row spectral sequence can be rewritten as the LES
\[
\cdots\to H_{p-1}K\to SH_pW\to H_p DW\to H_{p-2}K\to \cdots\]
Since nontrivial maps act by 0 on $K\iso\ker(W\to SW)$, the differentials vanish in the~complex
$K\tensor \epsilon_\ast$ computing $H_\ast K$, so $H_p K\iso K\tensor\epsilon_p$ 
(which was denoted $\widetilde{S}_{-p}K$ in [CEFN]). Therefore this LES becomes
\begin{equation}
\label{LES}
\cdots\to K\tensor \epsilon_{p-1}\to SH_pW\to H_p DW\to K\tensor \epsilon_{p-2}\to \cdots\end{equation}

Write $f_p(W)\coloneq \deg H_pW-p$.
\begin{maintheorem}[Church--Ellenberg]
\label{mainCE}
$f_p(W)\leq f_0(W)+f_1(W)$.
\end{maintheorem}
Our argument follows the proof in Li \cite{Li}, which depends on Li--Yu \cite{LiYu} (with ideas from Ramos \cite{Ramos}, Nagpal \cite{Nagpal}, and ultimately Church--Ellenberg \cite{CE}), but we remove the dependence on those earlier results; in particular, the characterizations of filtered FI-modules obtained in those papers are not necessary. Instead, we deduce \autoref{mainCE} from the long exact sequence \eqref{LES}.
\begin{proof}
We argue by induction on $f_0(W)$. Set $N\coloneq f_0(W)+f_1(W)$. In low degrees \eqref{LES} becomes
\begin{equation}
\label{Hone}
\cdots\to H_2 DW\to K\to SH_1 W\to H_1 DW\to 0\end{equation}
\begin{equation}
\label{Hzero}
0\to SH_0 W\to H_0DW\to 0\end{equation}From \eqref{Hone} and \eqref{Hzero} we see that $f_1(DW)\leq f_1(W)-1$ and $f_0(DW)=f_0(W)-1$. Applying \autoref{mainCE} to $DW$ by induction shows $f_p(DW)\leq N-2$. In particular, $\deg H_2DW\leq N$. Since $\deg SH_1W=f_1(W)\leq N$, from \eqref{Hone} we conclude $\deg K\leq N$, which implies $\deg K\otimes \epsilon_p\leq N+p$. Returning to \eqref{LES} for $p\geq 2$, we see the terms on either side of $SH_pW$ have degree bounded by $N+(p-1)$ and $N+p-2$ respectively, so $\deg SH_pW\leq N-1+p$. In other words, $f_p(W)\leq N$ as desired.
\end{proof}

\end{document}